\documentclass{article}

\usepackage{hyperref}
\usepackage{amsmath}
\usepackage{amsfonts}
\usepackage{amssymb}

\usepackage{graphicx}

\newtheorem{theorem}{Theorem}
\newtheorem{acknowledgement}[theorem]{Acknowledgement}

\begin{document}

\title{PVC Polyhedra}
\author{David Glickenstein\thanks{Partially funded by NSF DMS 0748283.}\\University of Arizona, Tucson, AZ 85721}
\maketitle

\begin{abstract}
We describe how to construct a dodecahedron, tetrahedron, cube, and octahedron
out of pvc pipes using standard fittings. 

\end{abstract}

\section{Introduction}

I wanted to build a huge dodecahedron for a museum exhibit. What better way to
draw interest than a huge structure that you can walk around and through? The
question was how to fabricate such an object? The dodecahedron is a fairly
simple solid, made up of edges meeting in threes at certain angles that form
pentagon faces. We could have them 3D printed. But could we do this with
``off-the-shelf''\ parts?

\begin{figure}
[ptb]
\begin{center}
\includegraphics
{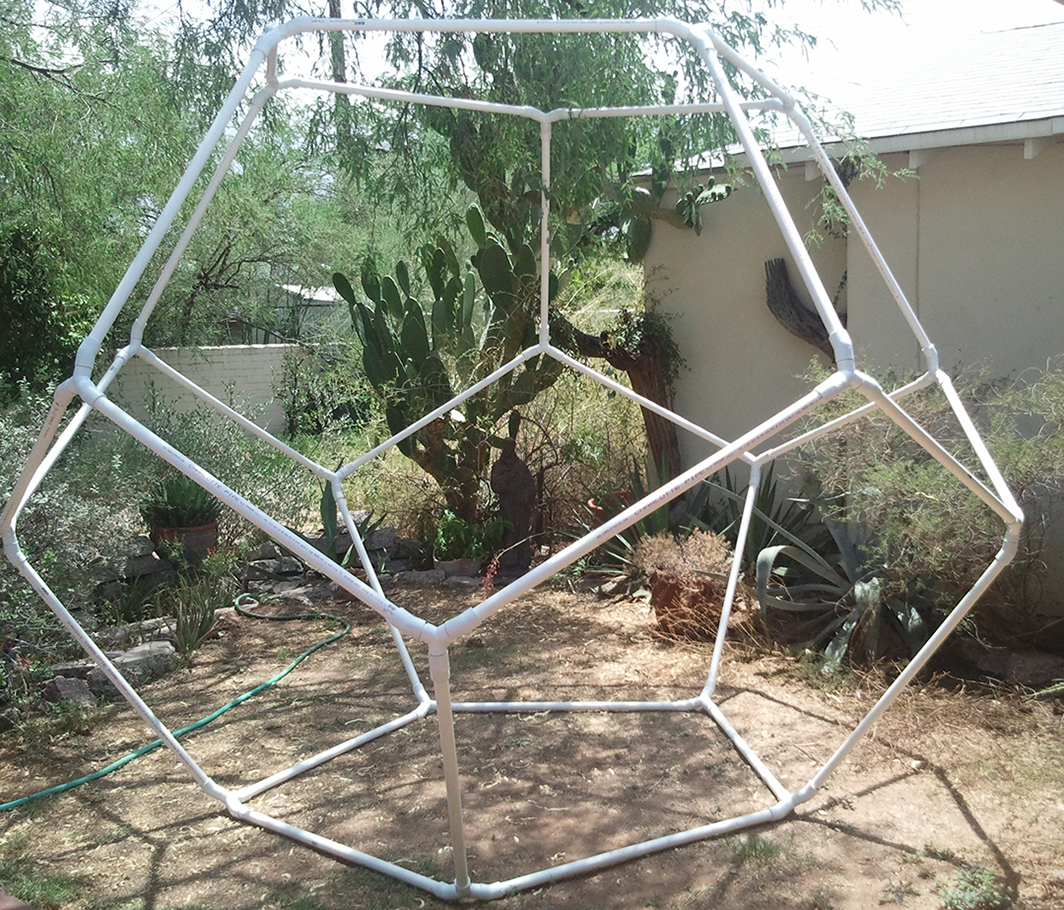}%
\caption{The PVC Dodecahedron. Picture courtesy of Bruce Bayly.}%
\label{fig:dodec}%
\end{center}
\end{figure}

The answer is yes, as seen in Figure \ref{fig:dodec}. The key fact to consider
is that each vertex corner consists of three edges coming together in a fairly
symmetric way. Therefore, we can take a connector with three pipe inputs and
make the corner a graph over it. In particular, if we take a connector that
takes three pipes each at 120 degree angles from the others (this is called a
``true wye'') and we take elbows of the
appropriate angle, we can make the edges come together below the center at
exactly the correct angles.%

\section{What is the correct elbow angle?}

\begin{figure}
[ptb]
\begin{center}
\includegraphics
{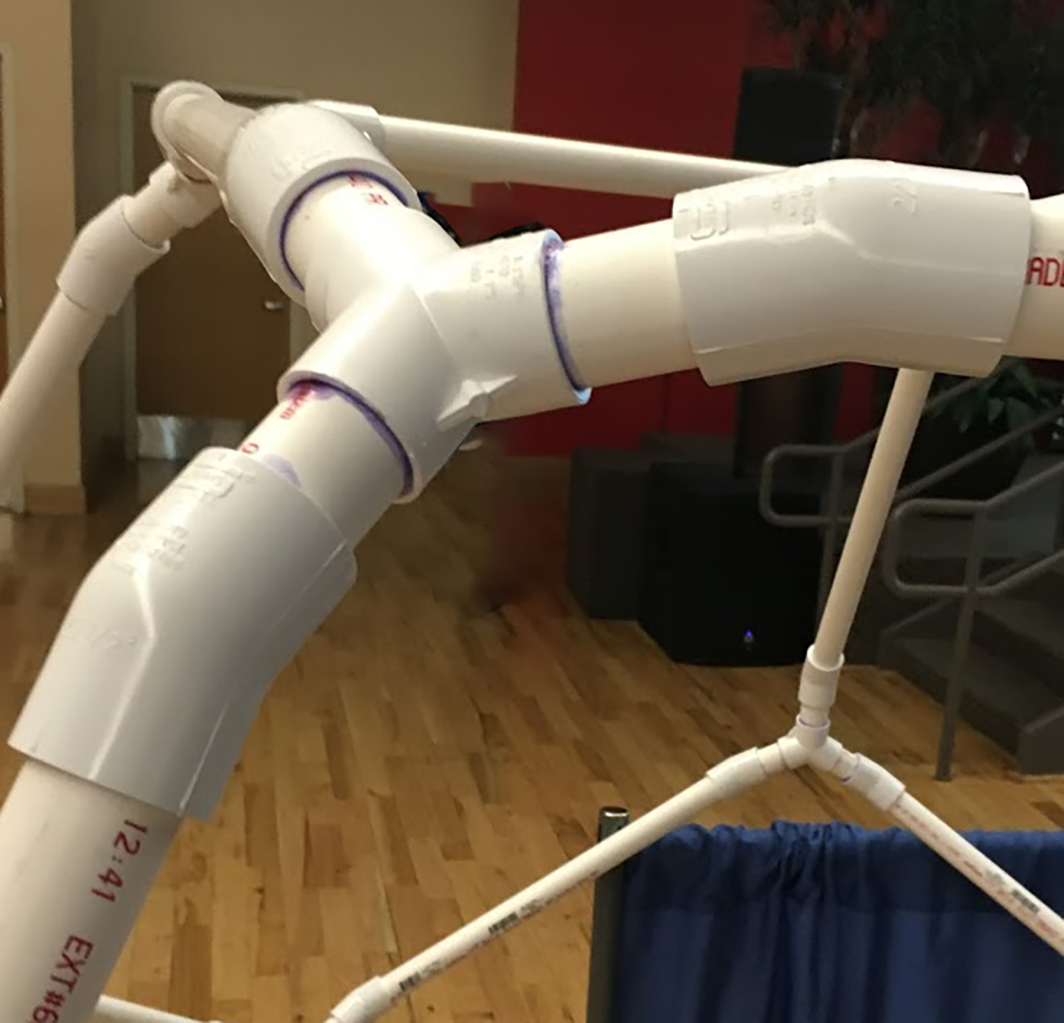}%
\caption{A close-up of the vertex joint for the dodecahedron.}%
\label{fig:dodecvertex}%
\end{center}
\end{figure}

Suppose the wye has three length 1 pipes connected to it and the actual vertex
$v$ is below the vertex $v_{W}$ of the wye. Let $X$ and $Y$ be the endpoints
of two of the wye segments. Then the segment $XY$ is of length $\sqrt{3}$
since it is opposite a 120 degree angle in the triangle $XYv_{W}.$ We know
that the angle $XvY$ at vertex $v$ has the angle of a regular pentagon, which
is $108$ degrees. Thus the length of segment $Xv$ is
\[
\frac{\sqrt{3}}{2}/\sin\left(  108%
{{}^\circ}%
\ast\frac{1}{2}\right)  =\frac{1}{2}\sqrt{15}-\frac{1}{2}\sqrt{3}=1.\,071.
\]
The needed angle at the elbow (in degrees) is
\[
\arccos\left(  \frac{1}{1.\,071}\right)  =20.\,98%
{{}^\circ}%
\]
We can get an elbow that is 22.5 degrees (which is half of 45 degrees). That
would be off by about 1.5 degrees, so the angles for the pentagons are
actually off by 3 degrees, which is $3/108=2.8\%.$ The length of $vv_{W}$ is
then
\[
\sqrt{\left(  \frac{1}{2}\sqrt{15}-\frac{1}{2}\sqrt{3}\right)  ^{2}-1}%
=\frac{3}{2}-\frac{1}{2}\sqrt{5}=0.382.
\]
Thus the actual vertices are $0.382$ times the distance from the elbow to the
center of the wye. Ours is 2.5 inches, and so the actual vertices are
$0.382\ast2.5=0.955$ inches from the center of the wye. Furthermore, given
whatever length we have of the tube, the actual vertices will be an additional
$1.071\ast2.5=2.\,678$ inches on each side. A close-up of the vertex joint can
be seen in Figure \ref{fig:dodecvertex}.%

\section{More platonic solids}

There are actually four platonic solids: the tetrahedron, the cube, the
octahedron, the dodecahedron, and the icosahedron. The vertices of a cube can
be purchased at a plumbing supply shop. We can also make the octahedron and a
tetrahedron with careful consideration.

For a regular tetrahedron, we compute similarly%
\[
\left\vert Xv\right\vert =\frac{\sqrt{3}}{2}/\sin\left(  30%
{{}^\circ}%
\right)  =\sqrt{\frac{\sqrt{6}}{2}}3
\]
and the needed elbow is $\arccos\left(  \frac{1}{\sqrt{3}}\right)  =54.7%
{{}^\circ}%
.$ If we use a three-way cube coupling, the angle with the vertex is
\[
\arccos\left(  \frac{2}{\sqrt{6}}\right)  =35.3%
{{}^\circ}%
.
\]
So with this we just need an elbow of
\[
54.7-35.3=19.4%
{{}^\circ}%
.
\]
This is not too far from 22.5 degrees.

\begin{figure}
[ptb]
\begin{center}
\includegraphics
{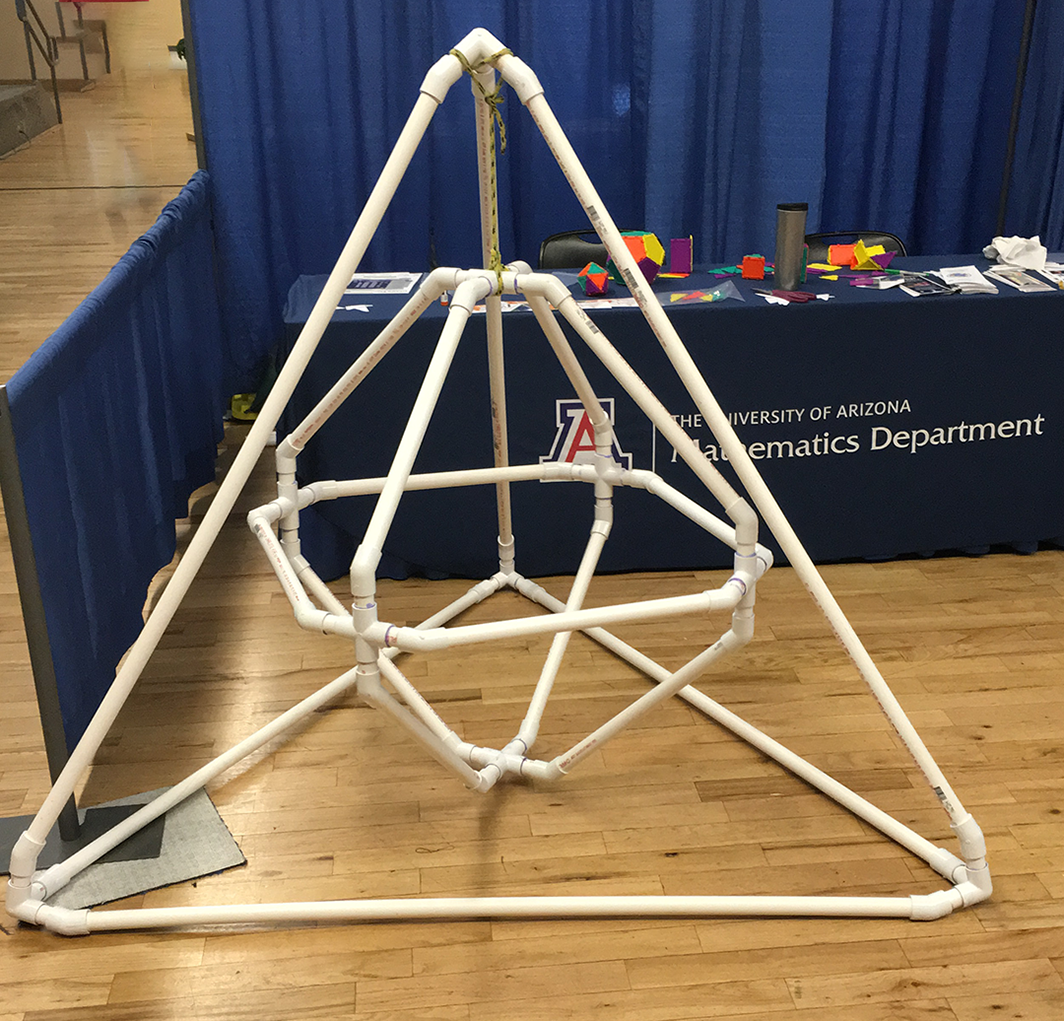}%
\caption{A PVC octahedron inside a PVC tetrahedron.}%
\label{fig:tetraocta}%
\end{center}
\end{figure}

For an octahedron, it can be verified that what is needed is simply a four-way
plus connector and 45 degree elbows. These can actually be purchased at most
large hardware/home improvement stores. The tetrahedron and octahedron can be
seen in Figure \ref{fig:tetraocta}%

\section{Make your own dodecahedron}

If you want to make your own dodecahedron, here is a quick inventory and size
calculator. Wikipedia is a good resource for some of these numbers \cite{wik}.
Since the dodecahedron has 20 vertices, you will need 20 true wyes and 60
twenty-two and half degree elbows. You will then need 60 small lengths of PVC
to connect the wyes and elbows, and you will want to use PVC glue to make the
connections strong. Then you will need 30 edges of roughly the same length (it
is very forgiving). To calculate the lengths, you need to decide how large you
want your dodecahedron to be.

Consider a dodecahedron with edge length $1.$ We have the following:

\begin{center}%
\begin{tabular}
[c]{|c|c|}\hline
inradius & $\frac{1}{20}\sqrt{250+110\sqrt{5}}\approx1.114$\\\hline
circumradius & $\frac{1}{4}\left(  \sqrt{15}+\sqrt{3}\right)  \approx
1.401$\\\hline
inradius of face & $\frac{1}{10}\sqrt{25+10\sqrt{5}}\approx0.688$\\\hline
circumradius of face & $\frac{1}{10}\sqrt{50+10\sqrt{5}}\approx0.851$\\\hline
dihedral angle & $\arccos\left(  -\frac{1}{\sqrt{5}}\right)  \approx
2.\,034\approx116.\,56^{\circ}$\\\hline
height of second row of vertices & $\frac{1}{2}\sqrt{5+2\sqrt{5}}\sin\left(
\pi-\arccos\left(  \frac{-1}{\sqrt{5}}\right)  \right)  \approx1.\,376$%
\\\hline
\end{tabular}

\end{center}

\noindent So to have a dodecahedron that sits $6$ ft. tall, we need edge
length of $6/2.228=2.\,693$ ft. Also, the height of the second row of vertices
(highest height to enter from) would be about $\left(  1.\,376\right)  \left(
2.693\right)  =3.\,706$ ft. Inside the dodecahedron we would have around
$\left(  0.688\right)  \left(  2.693\right)  \left(  2\right)  =3.\,706$ feet
in diameter to stand. The dodecahedron itself takes up a diameter of $\left(
1.401\right)  \left(  2.693\right)  \left(  2\right)  =7.\,546$ ft. The base
of the dodecahedron takes up a diameter of $\left(  0.851\right)  \left(
2.693\right)  \left(  2\right)  =4.\,584$ ft.

\section{Epilogue}

\begin{figure}
[ptb]
\begin{center}
\includegraphics
{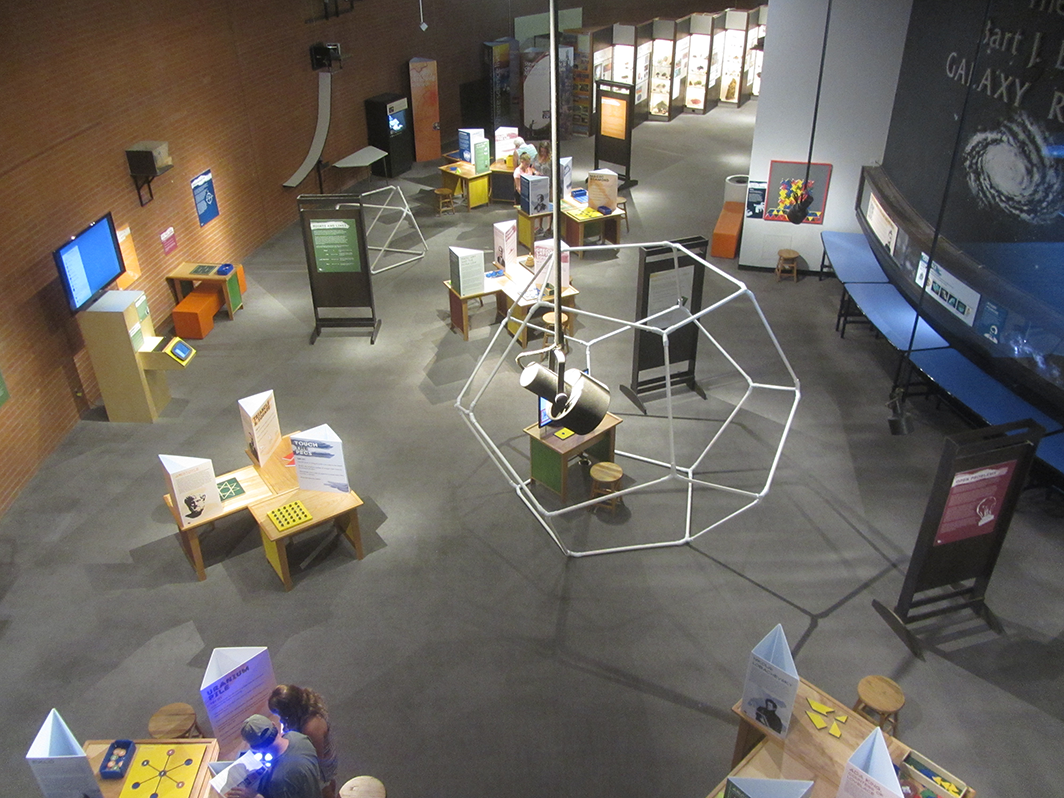}%
\caption{A birds-eye view of the exhibit Proofs, Puzzles, and Patterns:
Explore the World of Mathematics at Flandrau Science Center in Tucson, AZ.
Picture courtesy of Flandrau Science Center.}%
\label{fig:ppp}%
\end{center}
\end{figure}

The exhibit Proofs, Puzzles, and Patterns: Explore the World of Mathematics
\cite{PPP, Peiffer} opened at Flandrau Science Center in October 2015, as seen
in Figure \ref{fig:ppp}. The PVC dodecahedron was a hit, but in the end it was
not sturdy enough for the exhibit floor. The dodecahedron and other polyhedra
have found use in traveling exhibits by Bruce Bayly and his Mathematics Road
Show \cite{AMRS} as well as at some other enrichment locations. It is a great
group activity to put together the dodecahedron, and takes only around 10
minutes.%

\begin{acknowledgement}
The author would like to thank the staff at Flandrau Science Center for their
help and support in the construction of the Proofs, Puzzles, and Patterns
exhibit, especially Bill Plant, Shiloe Fontes, and Kellee Campbell, Neil
McSweeney, and Shipherd Reed; Greg McNamee for all of his work on the panels
for the exhibit; Marta Civil for her providing most of the puzzles aspect of
the exhibit; and Bruce Bayly, who is credited with the first picture and who
first assembled the PVC dodecahedron for the author in his own back yard, and
who has been a huge encouragement throughout.
\end{acknowledgement}

\end{document}